\baselineskip=14pt

\hsize=118mm 
\hoffset=20mm
\vsize=215mm
\voffset=15mm
\font\tencsc=cmcsc10

\def\ra{\rightarrow}

\def\fl{\forall}

\def\ot{\otimes}

\def\part{\partial}

\def\wt{\widetilde}

\def\d{\delta}
\def\g{\gamma}

\def\s{\sigma}

\def\ve{\varepsilon}
\def\vp{\varphi}

\def\D{\Delta}

\def\Lb{\Lambda}

\font\tenbb=msbm10
\font\sevenbb=msbm7
\font\fivebb=msbm5
\newfam\bbfam
\textfont\bbfam=\tenbb \scriptfont\bbfam=\sevenbb
\scriptscriptfont\bbfam=\fivebb
\def\bb{\fam\bbfam}

\def\Cb{{\bb C}}

\def\Zb{{\bb Z}}

\def\Hc{{\cal H}}

\def\Uc{{\cal U}}

\def\build#1_#2^#3{\mathrel{
\mathop{\kern 0pt#1}\limits_{#2}^{#3}}}

\def\hfl#1#2{\smash{\mathop{\hbox to 6mm{\rightarrowfill}}
\limits^{\scriptstyle#1}_{\scriptstyle#2}}}

\def\hfll#1#2{\smash{\mathop{\hbox to 6mm{\leftarrowfill}}
\limits^{\scriptstyle#1}_{\scriptstyle#2}}}

\def\boxit#1#2{\setbox1=\hbox{\kern#1{#2}\kern#1}%
\dimen1=\ht1 \advance\dimen1 by #1 \dimen2=\dp1 \advance\dimen2 by #1
\setbox1=\hbox{\vrule height\dimen1 depth\dimen2\box1\vrule}%
\setbox1=\vbox{\hrule\box1\hrule}%
\advance\dimen1 by .4pt \ht1=\dimen1
\advance\dimen2 by .4pt \dp1=\dimen2 \box1\relax}

\catcode`\@=11
\def\displaylinesno #1{\displ@y\halign{
\hbox to\displaywidth{$\@lign\hfil\displaystyle##\hfil$}&
\llap{$##$}\crcr#1\crcr}}

\def\ldisplaylinesno #1{\displ@y\halign{
\hbox to\displaywidth{$\@lign\hfil\displaystyle##\hfil$}&
\kern-\displaywidth\rlap{$##$}
\tabskip\displaywidth\crcr#1\crcr}}
\catcode`\@=12

\baselineskip=14pt

\hsize=118mm 
\hoffset=20mm
\vsize=215mm
\voffset=15mm

\def\Hc{\cal H}

\def\Uc{\cal U}

\font\tenbb=msbm10
\font\sevenbb=msbm7
\font\fivebb=msbm5
\newfam\bbfam
\textfont\bbfam=\tenbb \scriptfont\bbfam=\sevenbb
\scriptscriptfont\bbfam=\fivebb
\def\bb{\fam\bbfam}

\def\Cb{{\bb C}}

\def\Zb{{\bb Z}}

\def\d{\delta}
\def\D{\Delta}
\def\g{\gamma}

\def\L{\Lambda}

\def\s{\sigma}

\def\ve{\varepsilon}
\def\vp{\varphi}

\def\part{\partial}

\def\fl{\forall}

\def\ot{\otimes}

\def\wt{\widetilde}

\def\ra{\rightarrow}

\def\build#1_#2^#3{\mathrel{
\mathop{\kern 0pt#1}\limits_{#2}^{#3}}}

\def\xx{\vrule height 0.3em depth 0.2em width 0.3em}

\font\tenbb=msbm10
\font\sevenbb=msbm7
\font\fivebb=msbm5
\newfam\bbfam
\textfont\bbfam=\tenbb \scriptfont\bbfam=\sevenbb
\scriptscriptfont\bbfam=\fivebb
\def\bb{\fam\bbfam}

\def\Cb{{\bb C}}

\def\Zb{{\bb Z}}

\def\Hc{{\cal H}}

\def\Uc{{\cal U}}

\def\d{\delta}

\def\g{\gamma}

\def\s{\sigma}

\def\ve{\varepsilon}
\def\vp{\varphi}

\def\D{\Delta}

\def\Lb{\Lambda}

\def\fl{\forall}

\def\ot{\otimes}
\def\part{\partial}

\def\ra{\rightarrow}

\def\text{\hbox}

\def\build#1_#2^#3{\mathrel{
\mathop{\kern 0pt#1}\limits_{#2}^{#3}}}

\vglue 2cm

\font\tencsc=cmcsc10
\def\ra{\rightarrow}
\def\ot{\otimes}
\def\ve{\varepsilon}
\def\d{\delta}
\def\D{\Delta}
\def\ot{\otimes}
\def\wt{\widetilde}
\def\fl{\forall}
\def\s{\sigma}
\def\Lb{\Lambda}

\font\tenbb=msbm10
\font\sevenbb=msbm7
\font\fivebb=msbm5
\newfam\bbfam
\textfont\bbfam=\tenbb \scriptfont\bbfam=\sevenbb
\scriptscriptfont\bbfam=\fivebb
\def\bb{\fam\bbfam}

\def\Cb{{\bb C}}

\def\Hc{{\cal H}}

\def\boxit#1#2{\setbox1=\hbox{\kern#1{#2}\kern#1}%
\dimen1=\ht1 \advance\dimen1 by #1 \dimen2=\dp1 \advance\dimen2 by #1
\setbox1=\hbox{\vrule height\dimen1 depth\dimen2\box1\vrule}%
\setbox1=\vbox{\hrule\box1\hrule}%
\advance\dimen1 by .4pt \ht1=\dimen1
\advance\dimen2 by .4pt \dp1=\dimen2 \box1\relax}

\def\xx{\vrule height 0.5em depth 0.2em width 0.5em}


\centerline{\tencsc Cyclic cohomology and Hopf algebras}

\vglue 1cm

\centerline{A. Connes and H. Moscovici \footnote{(*)}{This material
 is based upon work supported by the National Science Foundation under Award No. DMS-9706886.}}
\bigskip \centerline{\it Dedicated to the memory of Mosh\'e Flato}
\vglue 1cm

\noindent{\bf Abstract}
\smallskip
We show by a direct computation that, for 
any Hopf algebra with a modulus-like character, the formulas first introduced 
in [CM] in the context of characteristic classes for actions of Hopf 
algebras, do define a cyclic module. This provides a natural generalization of 
Lie algebra cohomology to the general framework of Noncommutative Geometry, 
which covers the case of the Hopf algebra associated to 
$n$-dimensional transverse geometry [CM] as well as the function algebras
of the classical quantum groups.
\vglue 1cm

\noindent{\bf Introduction}
\medskip
We shall concentrate in this paper on the interplay between two basic concepts of 
Noncommutative geometry. The first is cyclic cohomology which plays the 
same role in Noncommutative geometry as De Rham cohomology plays in differential geometry.
 The second is Hopf algebras whose actions on noncommutative algebras are analoguous
 to Lie group actions on ordinary manifolds.
\smallskip
We shall show by a direct and elementary computation that, for 
any Hopf algebra with a modulus-like character, the natural 
cosimplicial module associated to the subjacent coalgebra structure 
can be upgraded to a {\it cyclic module} (or, rather, a module over
the {\it cyclic category} $\Lambda$, cf. [C, III.A]), by invoking both 
the product and the antipode. This cyclic module was first introduced 
in [CM] in the context of characteristic classes for actions of Hopf 
algebras, under a certain condition of existence of sufficiently 
nondegenerate actions, which made the verification of the axioms 
tautological. The fact that the latter assumption was superfluous has 
also been remarked by M.~Crainic [Cr], who recasted our construction in 
the framework of the Cuntz-Quillen formalism [CQ].
\vglue 1cm

\noindent{\bf I Cyclic cohomology and the cyclic category}

\medskip

The role of cyclic cohomology in Noncommutative geometry can be understood at several levels.
In its simplest guise it is a construction of invariants of K-theory, extending to the general framework the 
Chern-Weil characteristic classes of vector bundles and allowing for concrete computations on Noncommutative spaces.
For starters one should prove for oneself the following simple algebraic statement which extends to higher dimension
the obvious properties of the K-theory invariant provided by a trace $\tau$ on a noncommutative algebra $A$, by means of the equality,
$$
< E \, ,\tau >  = \tau(E) \leqno (1)
$$
for any idempotent $E, \, E^2 =E$, $E \in M_q(A)$, where the trace $\tau$ is extended to the algebra $M_q(A)$ of matrices over $A$ by,
$$
\tau ((a_{i,j})) \,  = \sum \tau(a_{i,i}) \leqno (2)
$$
To pass from this 0-dimensional situation to, say, dimension 2, one considers a trilinear form $\tau$
on the algebra $A$ which possesses the following compatibility with the algebra structure, reminiscent of the properties of a trace,
 
$$
\eqalign{\tau (a^1, a^2, a^0) \,  =  \tau(a^0, a^1, a^2) \quad \quad \quad \quad \quad \quad \quad \quad \quad \quad \quad  \cr
\tau (a^0a^1, a^2, a^3) \,-\tau (a^0,a^1a^2, a^3) \,+\tau (a^0, a^1, a^2a^3) \, -  \tau(a^3a^0, a^1, a^2)
 = 0
. \cr
\forall a^j \in A}. \leqno (3)
$$

\smallskip The statement then asserts that for each idempotent $E, \, E^2 =E$, $E \in M_q(A)$, the scalar
$\tau (E, E, E)$ remains constant when 
$E$ is deformed among the idempotents of $M_q(A)$.
This homotopy invariance of the resulting pairing between cyclic cocycles of arbitrary dimension (i.e. multilinear forms
on $A$ fulfilling the n-dimensional analogue of (3) ) and K-theory is the starting point of cyclic cohomology.
\smallskip At the conceptual level, cyclic cohomology is obtained as an Ext functor by linearisation of the non-additive
category of algebras and algebra homomorphisms ([C2]) using the additive category of $\L$-modules 
where $\L$ is the cyclic category. 
\smallskip The cyclic category is a small category which can be defined by generators and relations.
It has the same objects as the small category $\D$ of totally ordered finite sets and increasing maps.
Let us recall the presentation of $\D$. It has one object $[n]$ for each integer $n$, and is generated by 
faces $\d_i, [n-1] \ra [n]$ (the injection that misses $i$), and degeneracies $\s_j,[n+1] \ra [n] $ (the surjection
 which identifies $j$ with $j+1$), with the relations,
$$
\d_j \, \d_i = \d_i \, \d_{j-1} \ \hbox{for} \ i < j \, , \ \s_j \, \s_i = 
\s_i \, \s_{j+1} \qquad i \leq j \leqno (4)
$$
$$
\s_j \, \d_i = \left\{ \matrix{
\d_i \, \s_{j-1} \hfill &i < j \hfill \cr
1_n \hfill &\hbox{if} \ i=j \ \hbox{or} \ i = j+1 \cr
\d_{i-1} \, \s_j \hfill &i > j+1 \, . \hfill \cr
} \right.
$$
To obtain $\L$ one adds for each $n$ a new morphism $\tau_n, [n] \ra [n]$ such that,
$$
\matrix{
\tau_n \, \d_i = \d_{i-1} \, \tau_{n-1} &1 \leq i \leq n , &\tau_n \, \d_0 = 
\d_n \hfill \cr
\cr
\tau_n \, \s_i = \s_{i-1} \, \tau_{n+1} &1 \leq i \leq n , &\tau_n \, \s_0 = 
\s_n \, \tau_{n+1}^2 \cr
\cr
\tau_n^{n+1} = 1_n \, . \hfill \cr
} \leqno (5)
$$

The small category $\L$ is in fact best obtained as a quotient of the following 
category $E \, \L$. The latter has one object $(\Zb , n)$ for each $n$ and the 
morphisms $f : (\Zb , n) \ra (\Zb , m)$ are non decreasing maps, $(n,m \geq 
1)$
$$
f : \Zb \ra \Zb \ , \ f(x+n) = f(x)+m \qquad \fl \, x \in \Zb \, . \leqno (6)
$$
One has $\L = (E \, \L) / \Zb$ for the obvious action of $\Zb$ by translation.
The original definition of $\L$ (cf. [C2]) used homotopy classes of non decreasing maps 
from $S^1$ to $S^1$ of degree 1, mapping $\Zb / n$ to $\Zb / m$ and is trivially equivalent to the above.

\bigskip
Given an algebra $A$ one obtains a module over the small category $\L$ by assigning 
to each integer $n \geq 0$ the vector space $C^n$ of $n+1$-linear forms $\vp (x^0 , \ldots , x^n)$ on $A$, 
 while the basic operations are given by
$$
\matrix{
(\d_i \, \vp) (x^0 , \ldots , x^n) &=& \vp (x^0 , \ldots , x^i \, x^{i+1} , 
\ldots , x^n) \hfill &i=0,1,\ldots , n-1 \hfill \cr
\cr
(\d_n \, \vp) (x^0 , \ldots , x^n) &=& \vp (x^n \, x^0 , x^1 , \ldots , 
x^{n-1}) \hfill \cr
\cr
(\s_j \, \vp) (x^0 , \ldots , x^n) &=& \vp (x^0 , \ldots , x^j , 1 , x^{j+1} 
, \ldots , x^n) \hfill &j=0,1,\ldots , n \hfill \cr
\cr
(\tau_n \, \vp) (x^0 , \ldots , x^n) &=& \vp (x^n , x^0 , \ldots , x^{n-1}) 
\, . \hfill \cr
} \leqno (7)
$$
In the first two lines $\d_i : C^{n-1} \ra C^n$. In the third line $\s_i : 
C^{n+1} \ra C^n$. Note that $(\s_n \, \vp) (x^0 , \ldots , x^n) = \vp (x^0 , 
\ldots , x^n ,1)$, $(\s_0 \, \vp) (x^0 , \ldots , x^n) = \vp (x^0 , 1 , x^1 , 
\ldots ,$ $x^n)$.
\smallskip These operations satisfy the relations (4) and (5). This shows that any algebra $A$ gives rise 
canonically to a $\L$-module and allows ([C2][L]) to interpret the cyclic cohomology groups $HC^n(A)$ as 
$Ext^n$ functors. All of the general properties of cyclic cohomology
such as the long exact sequence relating it to Hochschild cohomology are shared by Ext of general $\L$-
modules and can be attributed to the equality of the classifying space $B\L$ of the small category $\L$
with the classifying space $BS^1$ of the compact one-dimensional Lie group $S^1$.

\vglue 1cm
\noindent{\bf II Characteristic classes for actions of Hopf algebras}

\medskip
Hopf algebras arise very naturally from their actions on noncommutative algebras. Given an algebra $A$, an 
action of the Hopf algebra $\Hc$ on $A$ is given by a linear map,
$$
{\Hc} \ot A \ra A \ , \ h \ot a \ra h(a) 
$$
satisfying $h_1 (h_2 \, a) = (h_1 \, h_2) (a) \quad \fl \, h_i \in {\Hc}, \, a \in A$ and
$$
h(ab) = \sum \, h_{(1)} \, (a) \, h_{(2)} \, (b)  \qquad \fl \, a,b  \in A \, , \ h \in 
{\Hc} \, . \leqno (1)
$$
where the coproduct of $h$ is,
$$
\D(h)= \, \sum \, h_{(1)} \, \ot \, h_{(2)} \, \leqno (2)
$$
In concrete examples, the algebra $A$ appears first, together with linear maps $A \ra A$ satisfying
 a relation of the form (1) which dictates the Hopf algebra structure. We refer to [CM] for an 
application of this construction to the leaf space of foliations.
\medskip\smallskip
\noindent The theory of characteristic classes for actions of $\Hc$ extends 
the construction ([C3]) of cyclic cocycles from a Lie algebra of 
derivations of a $C^*$ algebra $A$, together with an {\it invariant trace} 
$\tau$ on $A$.
\smallskip
In order to cover the nonunimodular case which does appear in the simplest examples, we
 fix a character $\d$ of $\Hc$ which will play the role of the module of locally compact groups.
\smallskip
 We then introduce the twisted antipode,
$$
\wt S (y) = \sum \, \d (y_{(1)}) \, S (y_{(2)}) \ , \ y \in {\Hc} \, , \ \D \, 
y = \sum \, y_{(1)} \ot y_{(2)} \, . \leqno (3)
$$
One has $\wt S (y) = S (\s (y))$ where $\s$ is the automorphism $\s= \, (\d \ot 1) \circ \D : {\Hc} \ra {\Hc}$.
\medskip

\noindent {\bf Definition 1.}
{\it We shall say that a trace $\tau$ on $A$ is $\d$-invariant under the action of $\Hc$ iff the following 
holds,}
$$
\tau (h(a)b) = \tau (a \, \wt S (h)(b)) \qquad \fl \, a,b  \in A \, , \ h \in 
{\Hc} \, . 
$$
 We have shown in ([CM]) that the definition of the cyclic 
complex $HC^*_{\d} ({\Hc})$ is uniquely dictated in such a way that the following holds,

\medskip

\noindent {\bf Proposition 2.} ([CM]) {\it Let $\tau$ be a $\d$-invariant trace 
on $A$, then the following defines a canonical map from $HC^*_{\d} ({\Hc})$ to 
$HC^* (A)$,
$$
\matrix{
\g (h^1 \ot \ldots \ot h^n) \in C^n (A) \, , \ \g (h^1 \ot \ldots \ot h^n) 
(x^0 , \ldots , x^n) = \cr
\cr
\tau (x^0 \, h^1 (x^1) \ldots h^n (x^n)) \, . \cr
}
$$
}
In ([CM]) we needed to assume the existence of enough such actions of ${\Hc}$
in order to check that the formulas were actually defining a cyclic module.
We shall show below by a direct and elementary computation that, for 
any Hopf algebra with a modulus-like character $\d$ as above, the construction of [CM]
does yield a cyclic module.

\vglue 1cm
\noindent{\bf III The cyclic module of a Hopf algebra}

\medskip
In this section we shall associate a cyclic complex (in fact a $ \Lb $-module, where $ \Lb $ is the cyclic category), to any Hopf algebra together with a character $\d$ such that the $\d$-twisted antipode is an involution.  The resulting cyclic cohomology appears to be the natural candidate for the analogue of Lie algebra cohomology in the context of Hopf algebras, where both the Hochschild cohomology (also called Sweedler cohomology) or the transposed (also called Harrison cohomology) give too naive results. 
\smallskip

Let $\Hc$ be a Hopf algebra (over $\Cb$) with unit map $\eta : \Cb 
\ra \Hc$, counit $\ve : \Hc \ra \Cb$ and antipode $S : \Hc \ra \Hc$, 
$$S * I = I * S = \eta \, \ve.$$
We fix a character $\d : \Hc \ra \Cb$, which will play the role of the 
modular function of a locally compact group. With the usual coproduct 
notation
$$
\D h = \sum_{(h)} h_{(1)} \ot h_{(2)} \quad , \quad h \in \Hc \, ,
$$
we introduce the $\d$-twisted antipode
$$
\wt S (h) = \sum_{(h)} \d (h_{(1)}) \ S (h_{(2)}) \quad , \quad h \in 
\Hc \, . \leqno (1)
$$

The elementary properties of $S$ imply immediately that $\wt S$ is an 
algebra antihomomorphism
$$
\matrix{
&\wt S (h^1 \, h^2) = \wt S (h^2) \, \wt S (h^1) \quad , \quad \fl \, 
h^1 , h^2 \in \Hc \cr \cr
&\wt S (1) = 1 \, , \hfill \cr
} \leqno (2)
$$
a coalgebra twisted antimorphism
$$
\D \, \wt S (h) = \sum_{(h)} S (h_{(2)}) \ot \wt S (h_{(1)}) \quad , 
\quad \fl \, h \in \Hc \, ; \leqno (3)
$$
and also that it satisfies
$$
\ve \circ \wt S = \d \, . \leqno (4)
$$

By transposing the standard simplicial operators underlying the 
\break Hochschild homology complex of an algebra, one associates to $\Hc$, 
viewed only as a coalgebra, the following natural cosimplicial 
module: $\{ \Hc^{\ot n} \}_{n \geq 1}$, with face operators $\d_i : 
\Hc^{\ot n-1} \ra \Hc^{\ot n}$,
$$
\matrix{
&\d_0 (h^1 \ot \ldots \ot h^{n-1}) = 1 \ot h^1 \ot \ldots \ot h^{n-1} 
\, , \hfill \cr \cr
&\d_j (h^1 \ot \ldots \ot h^{n-1}) = h^1 \ot \ldots \ot \D h^j \ot 
\ldots \ot h^n \, , \ \fl \, 1 \leq j \leq n-1 \, , \hfill \cr \cr
&\d_n (h^1 \ot \ldots \ot h^{n-1}) = h^1 \ot \ldots \ot h^{n-1} \ot 1 
\hfill \cr
} \leqno (5)
$$
and degeneracy operators $\s_i : \Hc^{\ot n+1} \ra \Hc^{\ot n}$,
$$
\s_i (h^1 \ot \ldots \ot h^{n+1}) = h^1 \ot \ldots \ot \ve (h^{i+1}) 
\ot \ldots \ot h^{n+1} \ , \ 0 \leq i \leq n \, . \leqno (6)
$$

In [CM, \S~7] the remaining two essential features of a Hopf algebra 
--{\it product} and {\it antipode} -- are brought into play, to define the 
{\it cyclic operators} $\tau_n : \Hc^{\ot n} \ra \Hc^{\ot n}$,
$$
\tau_n (h^1 \ot \ldots \ot h^n) = (\D^{n-1} \, \wt S (h^1)) \cdot h^2 
\ot \ldots \ot h^n \ot 1 \, . \leqno (7)
$$

\bigskip

\noindent {\bf Theorem 3.} {\it Let $\Hc$ be a Hopf algebra endowed with a character 
$\d \in \Hc^*$ such that the corresponding twisted antipode $(1)$ is 
an involution:
$$
\wt{S}^2 = I \, . \leqno (8)
$$
Then $\Hc_{\d}^{\natural} = \{ \Hc^{\ot n} \}_{n \geq 1}$ equipped 
with the operators given by $(5)$--$(7)$ defines a module over the
cyclic category $\Lambda$.}

\bigskip

{\it Proof.} One has to check the relations
$$
\matrix{
&\tau_n \, \d_i = \d_{i-1} \, \tau_{n-1} \ , \ 1 \leq i \leq n \, , 
\hfill \cr \cr
&\tau_n \, \d_0 = \d_n \, , \hfill \cr
} \leqno (9)
$$
$$
\matrix{
&\tau_n \, \s_i = \s_{i-1} \, \tau_{n+1} \ , \ 1 \leq i \leq n \, , 
\hfill \cr \cr
&\tau_n \, \s_0 = \s_n \, \tau_{n+1}^2 \, , \hfill \cr
} \leqno (10)
$$
$$
\tau_n^{n+1} = I_n \, . \leqno (11)
$$

It is the latter which poses a technical challenge. To size it up, 
let us first look at the case $n=2$.

\medskip

In what follows we shall only use the basic properties of the 
product, the coproduct, the antipode and of the twisted 
antipode (cf. (1)--(4)). We shall also employ the usual notational 
conventions for the Hopf algebra calculus (cf. [S]).

To begin with,
$$
\eqalign{
\tau_2 (h^1 \ot h^2) = & \ \D \, \wt S (h^1) \cdot h^2 \ot 1 = \cr
= & \ \sum \wt S (h^1)_{(1)} \, h^2 \ot \wt S (h^1)_{(2)} \cr
= & \ \sum S (h_{(2)}^1) \, h^2 \ot \wt S (h_{(1)}^1) \, . \cr
}
$$
Its square is therefore:
$$
\eqalign{
\tau_2^2 (h^1 \ot h^2) = & \ \sum S (S (h_{(2)}^1)_{(2)} \, 
h_{(2)}^2) \, \wt S (h_{(1)}^1) \ot \wt S (S (h_{(2)}^1)_{(1)} \, 
h_{(1)}^2) \cr
= & \ \sum S (S (h_{(2)(1)}^1) \, h_{(2)}^2) \, \wt S (h_{(1)}^1) \ot 
\wt S (S (h_{(2)(2)}^1) \, h_{(1)}^2) \cr
= & \ \sum S (h_{(2)}^2) \, (S \circ S) \, (h_{(2)(1)}^1) \, \wt S 
(h_{(1)}^1) \ot \wt S (h_{(1)}^2) \, (\wt S \circ S) \, 
(h_{(2)(2)}^1) \cr
= & \ \sum S (h_{(2)}^2) \ \hbox{\boxit{6pt}{$S(S(h_{(1)(2)}^1)) \, 
\wt S (h_{(1)(1)}^1)$}} \ \ot \wt S (h_{(1)}^2) \, \wt S (S 
(h_{(2)}^1)) \, . \cr
}
$$
The term in the box is computed as follows. With $k = h_{(1)}^1$, one 
has
$$
\eqalign{
\sum \, S (S(k_{(2)})) \, \wt S (k_{(1)}) = & \ \sum S (S(k_{(2)})) 
\, \d (k_{(1)(1)}) \, S (k_{(1)(2)}) \cr
= & \ \sum S(S(k_{(2)(2)}) \, \d (k_{(1)}) \, S (k_{(2)(1)}) = \cr
= & \ \sum \d (k_{(1)}) \, S \left( \sum k_{(2)(1)} \, S (k_{(2)(2)}) 
\right) \cr
= & \ \sum \d (k_{(1)}) \, S (\ve (k_{(2)}) \, 1) = \cr
= & \ \sum \d (k_{(1)}) \, \ve (k_{(2)}) = \d \left( \sum k_{(1)} \, 
\ve (k_{(2)}) \right) \cr
= & \ \d (k) \, . \cr
}
$$
It follows that
$$
\eqalign{
\tau_2^2 (h^1 \ot h^2) = & \ \sum S (h_{(2)}^2) \, \underbrace{\d 
(h_{(1)}^1) \ot \wt S (h_{(1)}^2) \, \wt S (S (h_{(2)}^1))} \cr
= & \ \sum S (h_{(2)}^2) \ot \wt S (h_{(1)}^2) \, \wt S (\wt S (h^1)) = \cr
= & \ \sum S (h_{(2)}^2) \ot \wt S (h_{(1)}^2) \, h^1 \, , \cr
}
$$
where we have used first (1) then (8). Thus
$$
\eqalign{
\tau_2^2 (h^1 \ot h^2) = & \sum  S (h_{(2)}^2) \ot \wt S (h_{(1)}^2) 
\cdot 1 \ot h^1 \cr
= & \ \D \, \sum \wt S (h^2) \cdot 1 \ot h^1 \, . \cr
}
$$
In a similar fashion,
$$
\eqalign{
\tau_2^3 (h^1 \ot h^2) = & \ \sum S (S(h_{(2)}^2)_{(2)}) \, \wt S 
(h_{(1)}^2) \, h^1 \ot \wt S (S (h_{(2)}^2)_{(1)}) \cr
= & \ \sum S (S(h_{(2)(1)}^2)) \, \wt S (h_{(1)}^2) \, h^1 \ot \wt S 
(S (h_{(2)(2)}^2)) \cr
= & \ \sum S (S(h_{(1)(2)}^2)) \, \wt S (h_{(1)(1)}^2) \, h^1 \ot \wt 
S (S (h_{(2)}^2)) \cr
= & \ \sum \d (h_{(1)}^2) \, h^1 \ot \wt S (S (h_{(2)}^2)) = \cr
= & \ \sum h^1 \ot \wt{S}^2 (h^2) = h^1 \ot h^2 \, . \cr
}
$$
\medskip

We now pass to the general case. With the standard conventions of 
notation,
$$
\eqalign{
\tau_n (h^1 \ot h^2 \ot \ldots \ot h^n) = \D^{(n-1)} \, \wt S (h^1) 
\cdot h^2 \ot \ldots \ot h^n \ot 1 \cr
= \sum S (h_{(n)}^1) \, h^2 \ot S (h_{(n-1)}^1) \, h^3 \ot \ldots \ot 
S (h_{(2)}^1) \, h^n \ot \wt S (h_{(1)}^1) \, . \cr
}
$$
Upon iterating once
$$
\eqalign{
\tau_n^2 (h^1 \ot \ldots \ot h^n) = & \ \sum S (S(h_{(n)}^1)_{(n)} \, 
h_{(n)}^2 )) \, S (h_{(n-1)}^1) \, h^3 \ot \cr
\ot & \ S (S (h_{(n)}^1)_{(n-1)} \, h_{(n-1)}^2)) \, S (h_{(n-2)}^1) 
\, h^4 \ot \ldots \cr
\ldots \ot & \ S (S (h_{(n)}^1)_{(2)} \, h_{(2)}^2)) \, \wt S 
(h_{(1)}^1) \ot \wt S (S (h_{(n)}^1)_{(1)} \, h_{(1)}^2) \cr
= & \ \sum S (h_{(n)}^2) \, S (S (h_{(n)(1)}^1)) \, S (h_{(n-1)}^1) 
\, h^3 \ot \cr
\ot & \ S (h_{(n-1)}^2) \, S (S (h_{(n)(2)}^1)) \, S (h_{(n-2)}^1) \, 
h^4 \ot \ldots \cr
\ldots \ot & \ S (h_{(2)}^2) \, S (S (h_{(n)(n-1)}^1)) \, \wt S 
(h_{(1)}^1) \ot \cr
\ot & \ \wt S (h_{(1)}^2) \, \wt S (S(h_{(n)(n)}^1)) = \cr
= & \ \sum S (h_{(n)}^2) \, S (h_{(n-1)}^1 \, S (h_{(n)}^1)) \, h^3 
\ot \cr
\ot & \ S (h_{(n-1)}^2) \, S (h_{(n-2)}^1 \, S (h_{(n+1)}^1)) \, h^4 
\ot \ldots \cr
\ldots \ot & \ S (h_{(2)}^2) \, S (S (h_{(2n-2)}^1 )) \cdot \wt S 
(h_{(1)}^1) \ot \cr
\ot & \ \wt S (h_{(1)}^2) \, \wt S ( S (h_{(2n-1)}^1 )) \, . \cr
}
$$
We pause to note that
$$
\sum h_{(n-1)}^1 \, S (h_{(n)}^1) = \sum h_{(n-1)(1)}^1 \, S 
(h_{(n-1)(2)}^1)
$$
equals
$$
\ve \, (h_{(n-1)}^1) \, 1 \, ,
$$
after resetting the indexation. Next
$$
\sum \ve (h_{(n-1)}^1) \, h_{(n-2)}^1
$$
gives $h_{(n-2)}^1$ after another resetting. In turn
$$
\sum h_{(n-2)}^1 \, S (h_{(n-1)}^1) 
$$
equals
$$
\ve (h_{(n-2)}^1) \, 1 \, ,
$$
and the process continues.

\smallskip

In the last step,
$$
\eqalign{
& \ \sum S (h_{(n)}^2) \, h^3 \ot S (h_{(n-1)}^2) \, h^4 \ot \ldots 
\cr
\ldots \ot & \ S (h_{(2)}^2) \ \hbox{\boxit{6pt}{$S(S(h_{(2)}^1)) \, \d 
(h_{(1)(1)}^1) \, S (h_{(1)(2)}^1)$}} \ \ot \cr
\ot & \ \wt S (h_{(1)}^2) \, \wt S ( S (h_{(3)}^1)) \cr
= & \ \sum S (h_{(n)}^2) \, h^3 \ot S (h_{(n-1)}^2) \, h^4 \ot \ldots
\ot S 
(h_{(2)}^2) \ot \wt S (h_{(1)}^2) \, h^1 \cr
= & \ S (h_{(n)}^2) \ot  S (h_{(n-1)}^2) \ot \ldots \ot \wt S (h_{(1)}^2)
\cdot 
h^3 \ot h^4 \ot \ldots \ot 1 \ot h^1 \cr
= & \ \D^{(n-1)} \, \wt S (h^2) \cdot h^3 \ot h^4 \ot \ldots \ot 1 \ot 
h^1 \, , \cr
}
$$
with the boxed term simplified as before. 

By induction, one obtains 
for any $j = 1, \ldots , n+1$,
$$
\tau_n^j (h^1 \ot \ldots \ot h^n) = \D^{n-1} \, \wt S (h^j) \cdot h^{j+1} 
\ot \ldots \ot h^n \ot 1 \ot \ldots \ot h^{j-1} \, ,
$$
in particular
$$
\tau_n^{n+1} (h^1 \ot \ldots \ot h^n) = \D^{n-1} \, \wt S (1) \cdot h^1 \ot 
\ldots \ot h^n = h^1 \ot \ldots \ot h^n \, .
$$

The verification of the compatibility relations (9), (10) is 
straightforward. Indeed, starting with the compatibility
with the face operators, one has:
$$
\eqalign{
\tau_n \, \d_0 (1 \ot h^1 \ot \ldots \ot h^{n-1}) = & \ \tau_n (1 \ot 
h^1 \ot \ldots \ot h^{n-1}) = \cr
= & \ \D^{n-1} \, \wt S (1) \cdot h^1 \ot \ldots \ot h^{n-1} \ot 1 
\cr
= & \ h^1 \ot \ldots \ot h^{n-1} \ot 1 \cr
= & \ \d_n (h^1 \ot \ldots \ot h^{n-1}) \, , \cr
}
$$
then
$$
\eqalign{
\tau_n \, \d_1 (h^1 \ot \ldots \ot h^{n-1}) = & \ \tau_n \,( \D \, h^1 
\ot h^2 \ot \ldots \ot h^{n-1}) \cr
= & \ \sum \tau_n (h_{(1)}^1 \ot h_{(2)}^1 \ot h^2 \ot \ldots \ot 
h^{n-1}) \cr
= & \ \sum \D^{n-1} \, \wt S (h_{(1)}^1) \cdot h_{(2)}^1 \ot h^2 \ot 
\ldots \ot h^{n-1} \ot 1 = \cr
= & \ \sum S (h_{(1)(n)}^1) \, h_{(2)}^1 \ot S (h_{(1)(n-1)}^1) \, 
h^2 \ot \ldots \cr
& \ \ot S (h_{(1)(2)}^1) \, h^{n-1} \ot \wt S (h_{(1)(1)}^1) \cr
= & \ \sum \ve (h_{(n)}^1) \, 1 \ot S (h_{(n-1)}^1) \, h^2 \ot \ldots 
\cr
& \ \ot S (h_{(1)}^1) \, h^{n-1} \ot \wt S (h_{(1)}^1) \cr
= & \ 1 \ot S (h_{(n-1)}^1) \, h^2 \ot \ldots \ot S (h_{(1)}^1) \, 
h^{n-1} \ot \wt S (h_{(1)}^1) \cr
= & \ \d_0 \, \tau_{n-1} \, (h^1 \ot \ldots \ot h^{n-1}) \, , \cr
}
$$
and so forth.

Passing now to degeneracies,
$$
\eqalign{
& \ \tau_n \, \s_0 (h^1 \ot \ldots \ot h^{n+1}) = \ve (h^1) \, \tau_n 
(h^2 \ot \ldots \ot h^{n+1}) = \cr
= & \ \ve (h^1) \, S (h_{(n)}^2) \, h^3 \ot \ldots \ot S (h_{(2)}^2) 
\, h^{n+1} \ot \wt S (h_{(1)}^2) \, , \cr
}
$$
and on the other hand
$$
\eqalign{
& \ \s_n \, \tau_{n+1}^2 (h^1 \ot \ldots \ot h^{n+1}) = \cr
= & \ \s_n \left( \sum S (h_{(n+1)}^2) \, h^3 \ot \ldots \ot S 
(h_{(2)}^2) \ot \wt S (h_{(1)}^2) \, h^1 \right) \cr
= & \ \sum \ve (\wt S (h_{(1)}^2) \, h^1) \, S (h_{(n+1)}^2) \, h^3 
\ot \ldots \ot S (h_{(2)}^2) \cr
= & \ \ve (h^1) \sum \d (h_{(1)}^2) \, S (h_{(n+1)}^2) \, h^3 \ot 
\ldots \ot S (h_{(2)}^2) \cr
= & \ \ve (h^1) \, S (h_{(n)}^2) \, h^3 \ot \ldots \ot S (h_{(2)}^2) 
\, h^{n+1} \ot \wt S (h_{(1)}^2) \, . \cr
}
$$
In the next step
$$
\eqalign{
& \ \tau_n \, \s_1 (h^1 \ot \ldots \ot h^{n+1}) =  \ve (h^2) \, 
\tau_n (h^1 \ot h^3 \ot \ldots \ot h^{n+1}) \cr
= & \ \ve (h^2) \cdot \D^{n-1} \, \wt S (h^1) \cdot h^3 \ot \ldots 
\ot h^{n+1} \ot 1 \, , \cr
}
$$
while on the other hand
$$
\eqalign{
& \ \s_0 \, \tau_{n+1} (h^1 \ot \ldots \ot h^{n+1}) = \cr
& \ \s_0 (S (h_{(n+1)}^1) \, h^2 \ot \ldots \ot S (h_{(2)}^1) \, 
h^{n+1} \ot \wt S (h_{(1)}^1)) \cr
= & \ \ve (h^2) \cdot \ve (h_{(n+1)}^1) \cdot S (h_{(n)}^1) \, h^3 
\ot \ldots \ot S (h_{(2)}^1) \, h^{n+1} \ot \wt S (h_{(1)}^1) \cr
= & \ \ve (h^2) \cdot S (h_{(n-1)}^1) \, h^3 \ot \ldots \ot S 
(h_{(2)}^1) \, h^{n+1} \ot \wt S (h_{(1)}^1) \, , \cr
}
$$
and similarly for $i = 2, \ldots n$.~\xx

\bigskip

The cohomology of the $(b,B)$-bicomplex corresponding to the cyclic 
module $\Hc_{\d}^{\natural}$ is, by definition, the {\it cyclic cohomology} 
$H \, C_{\d}^* (\Hc)$ of $\Hc$ relative to the modular character 
$\d$.
\bigskip
\noindent When ${\Hc} = {\Uc} ({\bf G})$ is the envelopping algebra of a Lie 
algebra, there is a natural interpretation of the Lie algebra cohomology,
$$
H^* ({\bf G} , \Cb) = H^* ({\Uc} ({\bf G}) , \Cb) 
$$
where the right hand side is the Hochschild cohomology with coefficients in 
the ${\Uc} ({\bf G})$-bimodule $\Cb$ obtained using the augmentation. In general, 
given a Hopf algebra ${\Hc}$ one can dualise (this is the construction of the Harrison complex),
the construction of the Hochschild 
complex $C^n ({\Hc}^* , \Cb)$ where $\Cb$ is viewed as a bimodule on ${\Hc}^*$ 
using the augmentation, i.e. the counit of ${\Hc}^*$. This gives the above 
operations: ${\Hc}^{\ot (n-1)} \ra {\Hc}^{\ot n}$, defining a cosimplicial 
space.
\smallskip
 When applied to the Hopf algebra ${\Hc}$ of functions on an affine algebraic group,
  this dual-Hochschild or Harrison cohomology gives 
simply the vector space of invariant twisted forms and ignores the group cohomology.
The second assertion of the following proposition shows that the cyclic cohomology  
$H \, C_{\d}^* ({\Hc})$, gives a highly nontrivial answer, thanks precisely to the action of the $B$ operator ([CM]).
\medskip

\noindent {\bf Proposition 4.} ([CM]) 1) {\it The periodic cyclic cohomology $H \, C_{\d}^* 
({\Hc})$, for ${\Hc} = {\Uc} ({\bf G})$ the envelopping algebra of a Lie 
algebra ${\bf G}$ is isomorphic to the Lie algebra homology $H_* ({\bf G} , 
\Cb)$ where $\Cb$ is a ${\bf G}$-module using the module $\d$ of $G$.}

\noindent 2) {\it  Let ${\Hc} = {\Uc} ({\bf G})_*$ be the Hopf algebra of polynomials in the coordinates
on an affine simply connected nilpotent group $G$. The periodic cyclic cohomology $H \, C_{\d}^* ({\Hc})$, is isomorphic
 to the Lie algebra cohomology of ${\bf G}$ 
with trivial coefficients.}

\smallskip

\noindent We refer to [CM] for the proof which was also reproduced in [Cr].
\medskip 

Finally we should point out that the existence of a twisted antipode
$\wt S$ of square one is still a partial unimodularity condition on a Hopf 
algebra. 
It was crucial for the results of [CM] that this 
condition is actually fulfilled for the Hopf algebra associated to 
$n$-dimensional transverse geometry of foliations. Moreover, as we shall
see now, it is also fulfilled by the most popular quantum groups.

First, recall that if $\Hc$ is {\it quasi-triangular}
(also called  {\it braided} ),  with universal $R$-matrix $R$, 
then (see e.g. [K])
$$
S^2 (x) = u \, x \, u^{-1} \, \, , x \in \Hc 
$$
with
$$
u = \sum \, S ( R^{(2)}) \, R^{(1)} \ , \quad \varepsilon (u) = 1
$$
and
$$
\D \, u = (R_{21} \, R)^{-1} \, (u \otimes u) \, .
$$

\noindent If in addition $u \, S(u) = S(u) \, u$, which is central, 
has a central square root $\theta$, such that
$$
\Delta (\theta) = (R_{21} \, R)^{-1} \, (\theta \otimes \theta) \ , \quad 
\varepsilon (\theta) = 1 \ , \quad S(\theta) = \theta \, ,
$$
then $\Hc$ is called a {\it ribbon} algebra. 
Any braided Hopf algebra $\Hc$ can be canonically embedded in a
ribbon algebra (cf. [RT]):
$$
\Hc (\theta) = \Hc \, [\theta] / (\theta^2 - u \, S(u))
$$
If $\Hc$ is a ribbon algebra, then
$$
\delta = \theta \, u^{-1} \, ,
$$
is a group-like element:
$$
\Delta \, \delta = \delta \otimes \delta \ , \quad \varepsilon (\delta) = 1 \ , 
\quad S(\delta) = \delta^{-1} \, .
$$
Defining 
$$
\wt S = \delta \cdot S \, ,
$$
one obtains a twisted antipode which satisfies the property
${\wt S}^2 = 1 \, .$
\noindent Indeed,
$$
\matrix{
{\widetilde S}^2 (x) &= &\delta \, S (\delta \, S (x)) = \delta \, S^2 (x) \, 
\delta^{-1} \hfill \cr
&= &\delta \, u \, x \, u^{-1} \, \delta^{-1} = \theta \, x \, \theta^{-1} = x 
\, , \hfill \cr
}
$$
because $\theta$ is central.
\medskip

By dualizing the above definitions one obtains the notion of a
{\it cobraided}, resp.  {\it coribbon} algebra. Among the most
prominent examples of coribbon algebras are the function
algebras of the classical quantum groups  $GL_q (N)$, $SL_q (N)$, $SO_q (N)$, 
$O_q (N)$ and $Sp_q (N)$.
  
For a coribbon algebra $\Hc$, the analogue of
the above {\it ribbon group-like element} $\delta$  is
the {\it ribbon character} $\delta \in \Hc^*$. The corresponding 
twisted antipode 
$${\wt S} = \delta * S \, , $$
satisfies again the condition ${\widetilde S}^2 = 1$.

\vglue 1cm

\noindent {\bf References}

\bigskip

\item{[C1]} A.~Connes, {\it Noncommutative Geometry}, London-New York, 
Academic Press, 1994.
\medskip

\item{[C2]} A. Connes : Cohomologie cyclique et foncteur $Ext^n$. {\it  C.R. Acad. Sci. 
Paris}, Ser.I Math {\bf 296} (1983).
\medskip

\item{[C3]} A. Connes : $C^*$ alg\`ebres et g\'eom\'etrie differentielle. {\it  C.R. Acad. Sci. 
Paris}, Ser.~A-B {\bf 290} (1980).
\medskip

\item{[CM]} A.~Connes and H.~Moscovici, Hopf Algebras, Cyclic 
Cohomology and the Transverse Index Theorem, {\it Commun. Math. 
Phys.} {\bf 198} (1998), 199-246.
\medskip

\item{[Cr]} M.~Crainic, Cyclic Cohomology of Hopf Algebras and a 
Noncommutative Chern-Weil Theory, {\it Preprint} QA/9812113.
\medskip

\item{[CQ]} J.~Cuntz and D.~Quillen, Cyclic homology and singularity,
{\it J. Amer. Math. Soc.} {\bf 8} (1995), 373-442.
\medskip

\item{[K]} C.~Kassel, {\it Quantum Groups}, Springer-Verlag, 1995.
\medskip

\item{[L]} J.L.~Loday, {\it Cyclic Homology}, Springer, Berlin, Heidelberg, 
New York, 1998.
\medskip

\item{[RT]} N. Yu.~Reshetikhin and V. G.~Turaev, Ribbon graphs and their
invariants derived from quantum groups,
{\it Commun. Math. Phys.} {\bf 127} (1990), 1-26.
\medskip

\item{[S]} M.E.~Sweedler, {\it Hopf Algebras}, W.A.~Benjamin, Inc., 
New York, 1969.

\bye